\documentclass[graybox]{svmult}

\usepackage{mathptmx}       
\usepackage{helvet}         
\usepackage{courier}        
\usepackage{type1cm}        
\usepackage{makeidx}         
\usepackage{graphicx}        
\usepackage{multicol}        
\usepackage[bottom]{footmisc}
\usepackage[french]{babel} 

\usepackage{amsmath}
\usepackage{amssymb}

\begin{document}

\title*{Une version  du th\'eor\`eme d'Amer et Brumer pour les z\'ero-cycles}
\author{Jean-Louis Colliot-Th\'el\`ene et Marc Levine}
\institute{Jean-Louis Colliot-Th\'el\`ene \at C.N.R.S., Math\'ematiques, B\^atiment 425, Universit\'e Paris-Sud, F-91405 Orsay, France, \email{jlct@math.u-psud.fr}
\and Marc Levine \at Department of Mathematics, Northeastern University,
360 Huntington Avenue,
Boston, MA 02115, USA  \
 and
\
Universit\"at Duisburg-Essen Fakult\"at Mathematik, Campus Essen, D-45117 Essen, Germany, 
 \email{marc.levine@uni-due.de}}
%
%

 \motto{\`A Parimala}

\maketitle

  
\newcommand{\CH}{\operatorname{CH}}

\def\Bbb{\mathbb}
\def\lgr{\longrightarrow}
\def\ra{\rightarrow}
\def\N{{\Bbb N}}
\def\Z{{\Bbb Z}}
\def\Q{{\Bbb Q}}
\def\R{{\Bbb R}}
\def\C{{\Bbb C}}
\def\F{{\Bbb F}}
\def\br{{\rm Br}}
\def\k{{\overline k}}

\def\pgcd{{\rm pgcd}}
\def\pic{{\rm Pic}}
\def\dim{{\rm dim}}
\def\spec{{\rm Spec} \ }
\def\Gr{{\rm Gr}}

\def\A{{\Bbb A}}
 \def\P{{\Bbb P}}
 \def\oi{\hskip1mm {\buildrel \cong \over \rightarrow} \hskip1mm}

 \motto{A Parimala}

\section{Introduction}

Soit $k$ un corps.
Soient $f$ et $g$ deux formes {\it quadratiques} \`a coefficients dans $k$,
en $n+1$ variables. Soit $t$ une variable.
On sait que le syst\`eme de formes $f=g=0$ a un z\'ero non trivial sur $k$
si et seulement si la forme quadratique $f+tg$ sur le corps $k(t)$
a un z\'ero non trivial (M. Amer \cite{A}, A. Brumer \cite{B}, voir \cite[III, Prop. 17.14]{EKM}).

Dans cette note, dont une version primitive fut con\c cue \`a Boston en avril 1991, nous montrons
que si l'on consid\`ere les z\'ero-cycles de degr\'e 1 plut\^ot que
les points rationnels, il existe en {\it tout degr\'e} $d\geq 2$ une version
 de ce r\'esultat pour un syst\`eme de deux formes 
 de  m\^eme  degr\'e  $d$  (Th\'eor\`eme \ref{extheoreme1}), et  des versions pour un syst\`eme 
quelconque de formes de m\^eme degr\'e $d$ (Th\'eor\`emes \ref{extheoreme2}  et \ref{extheoreme3}).

Pour  $K/k$ une extension de corps et   $X$ un $k$-sch\'ema,
on note $X_{K}=X\times_{k}K$.
Une vari\'et\'e alg\'ebrique   sur un corps $k$ est un $k$-sch\'ema 
s\'epar\'e de type fini.

\section
{\bf Indice et indice r\'eduit}

\bigskip

\begin{definition}
Soient $k$ un corps et $X$ une $k$-vari\'et\'e alg\'ebrique.
\`A tout point ferm\'e $P \in X$, de corps r\'esiduel $k(P)$  on associe son degr\'e $[k(P):k] \in \N$.
 L'indice $I(X/k)$ de la $k$-vari\'et\'e $X$ est par d\'efinition le plus grand commun diviseur (pgcd)  des degr\'es  des points ferm\'es.
 C'est aussi le pgcd des degr\'es $[L:k]$ des extensions finies de corps $L/k$ avec $X(L)\neq \emptyset$.
 \end{definition}
 
  De fa\c con \'evidente, on a $I(X/k)=1$ si et seulement si $X$ 
  poss\`ede un z\'ero-cycle (cycle de dimension z\'ero)
 $\sum_{P}n_{P}P$ (avec $n_{P} \in \Z$) de degr\'e $\sum_{P}n_{P}[k(P):k]=1$.
 On appelle indice r\'eduit de $X/k$, et on note 
 $I_{red}(X/k)$ le produit des nombres premiers qui divisent $I(X/k)$.
 De fa\c con triviale, $I(X/k)=1$  si et seulement si $I_{red}(X/k)=1$.

 \begin{lemma}
 \label{exlemme0}
   Soient $k$ un corps et $X$ une $k$-vari\'et\'e. 
 \begin{enumerate}
 \renewcommand{\theenumi}{\alph{enumi}}
 \renewcommand{\labelenumi}{(\theenumi)}
  \item  Si $X=Y \cup Z$ est la r\'eunion ensembliste  de deux sous-$k$-vari\'et\'es,
 alors   $I(X/k)=\pgcd(I(Y/k),I(Z/k)).$ En particulier, l'indice de $X/k$
 est \'egal \`a l'indice de sa sous-$k$-vari\'et\'e r\'eduite.
 
\item Si $X=Y \cup Z$ est la r\'eunion ensembliste  de deux sous-$k$-vari\'et\'es,
 alors $I_{red}(X/k)=\pgcd(I_{red}(Y/k),I_{red}(Z/k)).$

\item Si $K=k(t_{1}, \cdots, t_{n})$ est une extension transcendante pure de $k$,
 alors $I(X/k)=I(X_{K}/K)$ et $I_{red}(X/k)=I_{red}(X_{K}/K)$.
 \end{enumerate}
 \end{lemma}
 \begin{proof}
  Seul le point (c)   requiert une explication.
 Soit $K=k(t)$. Si $L/k$ est une extension finie de corps avec
 $X(L)\neq \emptyset$ alors $X_{k(t)}(L(t)) \neq \emptyset$.
  Ainsi $I(X_{K}/K)$ divise  $I(X/k)$. 
 Soit $P$ un point ferm\'e de $X_{K}$ de degr\'e $n$.
 On a donc une extension finie de corps $L/K$ et une
  $K$-immersion ferm\'ee
 $\spec L \hookrightarrow X_{K}$. Le corps $L$ est
 le corps des fonctions d'une $k$-courbe normale $C$,
 finie sur $\spec k[t]$. Il existe un ouvert non vide $U \subset \spec k[t]$
 tel que la restriction $C_{U}/U$ soit finie de degr\'e $n$, et qu'il
 existe une $U$-immersion ferm\'ee $C_{U} \hookrightarrow X\times_{k}U$.
 Si le corps $k$ est infini, on choisit un $k$-point $R \in U(k)$.
 La fibre de $C_{U}/U$ au-dessus de ce $k$-point est 
 le spectre d'une $k$-alg\`ebre finie de degr\'e $n$ qui
 admet une $k$-immersion dans $X$.
 Une telle situation d\'efinit un z\'ero-cycle effectif de degr\'e $n$
 sur la $k$-vari\'et\'e $X$ (voir \cite[Appendices A1, A2, A3]{F}).
 Si le corps $k$ est fini, il existe un z\'ero-cycle 
 $\sum n_{i} R_{i}$ tel  que tous les points ferm\'es
 $R_{i}$ soient dans $U \subset \spec  k[t]$
 et que $\sum_{i} n_{i} [k(R_{i}):k]=1$.
 A chaque $R_{i}$ on associe par la m\'ethode ci-dessus  un z\'ero-cycle effectif  $z_{i}$
 de degr\'e
 $ n[k(R_{i}):k]$ sur la $k $-vari\'et\'e $X$.
 Le z\'ero-cycle $\sum_{i}n_{i}.z_{i}$
 est alors de degr\'e $n$ sur la $k$-vari\'et\'e $X$. Ainsi $I(X/k)$ divise $I(X_{K}/K)$.
 L'\'enonc\'e sur les indices r\'eduits r\'esulte imm\'ediatement
 de celui sur les indices.
\qed \end{proof}

 \section{Syst\`eme de deux formes}

\begin{theorem}
\label{extheoreme1}
Soient $k$ un corps, $f$ et $g$ deux formes de degr\'e $d\geq  1$ en $n+1\geq 3$
variables, non toutes deux nulles. Soient $t$ une variable et $K=k(t)$.

L'indice r\'eduit de la $K$-hypersurface $W \subset \P^n_{K}$
 d\'efinie par $f+tg=0$   co\"{\i}ncide avec l'indice r\'eduit de la $k$-vari\'et\'e $X \subset\P^n_{k}$ d\'efinie par $f=g=0$. 
 
 En particulier, la $K$-hypersurface $f+tg=0$ poss\`ede un z\'ero-cycle de degr\'e 1 si et seulement si la $k$-sous-vari\'et\'e $X$ de $\P^n_{k}$ d\'efinie par $f=g=0$ poss\`ede un z\'ero-cycle de degr\'e 1.
 \end{theorem}
\begin{proof}
Pour les $k$-vari\'et\'es, on omet l'indice $k$. 
Si  $L/k$ est une extension finie de corps  avec $X(L)\neq \emptyset$, alors $W(L(t)) \neq \emptyset$.
Ainsi    l'indice $I(W/K)$  divise $I(X/k)$, et donc l'indice r\'eduit $I_{red}(W/K)$ 
 divise  l'indice r\'eduit $I_{red}(X/k)$.
D'apr\`es le lemme \ref{exlemme0}~(c), pour \'etablir l'\'enonc\'e on peut remplacer le corps
$k$ par une extension transcendante pure. On supposera donc le corps $k$ infini.

\medskip

{\it Supposons d'abord que les formes
$f$ et $g$ n'ont pas de facteur commun non constant.}

Soit $V=\P^n \setminus X$.
Notons $I=I(X/k)$. De la suite exacte de localisation (\cite[I. Prop. 1.8]{F}) :
$$ \CH_{0}(X) \to \CH_{0}(\P^n) \to \CH_{0}(V) \to 0$$
et du fait que le degr\'e sur $k$ d\'efinit un isomorphisme $\CH_{0}(\P^n) \simeq \Z$,
on tire l'\'egalit\'e $\CH_{0}(V)=\Z/I$.
Soit $Z \subset \P^1 \times \P^n$ la $k$-vari\'et\'e d\'efinie par
$$ \lambda f+\mu g=0.$$
Via la  projection $\P^1 \times \P^n \to \P^n$, c'est l'\'eclat\'ee
de $\P^n$ le long de l'intersection compl\`ete $X \subset  \P^n$
([F], Appendix A, Remark A.6).
Soit $U \subset \P^1 \times \P^n$ l'ouvert compl\'ementaire de $Z$.
La projection $q: \P^1 \times \P^n \to \P^n$ induit un morphisme $q_{U}: U \to V$
qui fait de $U$ un fibr\'e en droites affines sur $V$.
On en d\'eduit un isomorphisme $q_{U}^* :  \Z/I=\CH_{0}(V) \oi \CH_{1}(U)$.
On a la suite exacte de localisation  (\cite[I. Prop. 1.8]{F}) :
$$ \CH_{1}(Z) \to \CH_{1}( \P^1 \times \P^n) \to \CH_{1}(U) \to 0.$$
Par ailleurs les poussettes associ\'ees \`a la projection $p :  \P^1\times \P^n \to \P^1$ et  \`a la projection $q: \P^1 \times \P^n \to \P^n$
induisent   un isomorphisme 
$$(p_{*},q_{*}) : \CH_{1}(\P^1 \times \P^n) \oi \CH_{1}(\P^1) \oplus \CH_{1}(\P^n) 
= \Z \oplus \Z.$$
Notons $$i : \CH_{1}(Z) \to \CH_{1}( \P^1 \times \P^n)  \oi \CH_{1}(\P^1) \oplus \CH_{1}(\P^n) 
= \Z \oplus \Z$$
l'application compos\'ee.
Sur $Z$, on trouve les 1-cycles suivants.

Le corps $k$ \'etant  infini, et les formes $f$ et $g$ sans facteur commun, on trouve un cycle $z_{1}=\P^1 \times R$, o\`u $R$ est un z\'ero-cycle effectif de degr\'e $d^2$ sur $X$ par intersection avec un espace lin\'eaire de codimension 2 convenable.
L'image par $i$ de $z_{1} $ est $(d^2,0)$.  

 L'hypoth\`ese que $f$ et $g$ n'ont pas de facteur commun assure que la
 $k$-vari\'et\'e $X \subset \P^n$ est de codimension~2.
 Soient $\Gr(1,\P^n)$ la grassmannienne des droites dans $\P^n$ et $E \subset \Gr(1,\P^n) \times \P^n$
  la vari\'et\'e d'incidence.
 Soient $r_{1},r_{2}$ les deux projections induites sur $E$.
 L'image r\'eciproque $r_{2} ^{-1}(X)$  est de codimension au moins 2
 dans $E$,  le morphisme $r_{1}$ a ses fibres de dimension 1, donc l'adh\'erence de 
 $r_{1}(r_{2} ^{-1}(X))$ dans  $\Gr(1,\P^n) $ est de codimension au moins 1.
 Le corps $k$ \'etant infini, et la $k$-vari\'et\'e $\Gr(1,\P^n) $ $k$-birationnelle
 \`a un espace projectif,
 on peut donc trouver
 une $k$-droite  $L=\P^1 \subset \P^n$  qui ne rencontre pas $X$.
  Soit $q_{Z} : Z \to \P^n$ la restriction de $q$ \`a $Z$.
 Le morphisme $q_{Z}$ induit un isomorphisme $Z\setminus q^{-1}(X) \oi  \P^n \setminus X$.
 Soit $z_{2} \subset Z$ l'image r\'eciproque de $L$ par cet isomorphisme. Ceci d\'efinit  un
 1-cycle sur $Z \subset \P^1 \times \P^n$, dont 
 l'image par $i$   est  $(d,1)$. En effet, ce 1-cycle est donn\'e par
 $L \to \P^1 \times \P^n$, o\`u la projection sur le second facteur est l'inclusion 
lin\'eaire $l : L \subset \P^n$, 
et o\`u la projection sur le premier facteur est
donn\'ee par $(-g(l),f(l))$ (on a not\'e  ici $l$ un syst\`eme de $n+1$ formes lin\'eaires en deux variables).
Comme $X$ ne rencontre pas~$l$, le couple $(f(l),g(l))$ de formes
homog\`enes de degr\'e $d$  n'a pas de z\'ero commun.

Soit $s=I(W/K)$.
L'hypersurface $f+tg=0$ sur le corps $K$ poss\`ede un z\'ero-cycle de degr\'e $s$.
L'adh\'erence d'un tel z\'ero-cycle dans $Z \subset \P^1 \times \P^n$ d\'efinit un 1-cycle $z_{3}$
 dont l'image par $i$ est de la forme $(s,a)$ pour un certain entier $a \in \Z$.
 
 Le quotient de $\Z \oplus \Z$ par le groupe engendr\'e par $(d^2,0), (d,1), (s,a)$ 
est annul\'e par  l'entier
 $\pgcd(d^2,s-ad)$.  
 De la suite de localisation on conclut  que, pour un certain
 entier $a \in \Z$, l'entier  $I=I(X/k)$ divise  $\pgcd(d^2,I(W/K)-ad)$. 
 Ainsi $I_{red}(X/k)$ divise $I(W/K)$.
 Comme par ailleurs $I_{red}(W/K)$ divise $I_{red}(X/k)$, on conclut  
 $$I_{red}(X/k)=I_{red}(W/K).$$

 \medskip

{\it Supposons maintenant que $f=h.f_{1}$ et $g=h.g_{1}$ avec $f_{1}$ et $g_{1}$ homog\`enes de m\^eme degr\'e  
sans facteur commun non constant et $h$ homog\`ene non constant. } 

Soit $X \subset \P^n_{k}$, resp. $X_{1} \subset \P^n_{k}$, resp. $X_{2} \subset \P^n_{k}$ la $k$-vari\'et\'e d\'efinie par $f=g=0$, resp. par $f_{1}=g_{1}=0$, resp. par $h=0$.
Soit $W\subset \P^n_{K}$, resp. $W_{1}/K$ la vari\'et\'e d\'efinie par $f+tg=0$, resp. $f_{1}+tg_{1}=0$.
On a
$$I_{red}(W/K)= \pgcd(I_{red}(W_{1}/K),I_{red}(X_{2,K}/K)) =  \pgcd(I_{red}(W_{1}/K),I_{red}(X_{2}/k)) $$
d'apr\`es le lemme \ref{exlemme0}
et
$$\pgcd(I_{red}(W_{1}/K),I_{red}(X_{2}/k)) = \pgcd(I_{red}(X_{1}/k),I_{red}(X_{2}/k)) =I_{red}(X/k),$$
la premi\`ere \'egalit\'e r\'esultant de $I_{red}(W_{1}/K)=I_{red}(X_{1}/k)$ \'etabli ci-dessus,
 la seconde \'egalit\'e   provenant du lemme \ref{exlemme0}.
Ceci ach\`eve la d\'emonstration. 
\qed \end{proof}

Le th\'eor\`eme \ref{extheoreme1}  se g\'en\'eralise \`a un nombre quelconque de formes. Il y a en fait deux g\'en\'eralisations.
Nous aurons besoin du lemme suivant.

\begin{lemma}\label{newlemma}
 Soient $k$ un corps, $f_{0},f_{1},\dots,f_{r}$ des formes  \`a coefficients dans $k$,  
 de degr\'e $d\geq  1$, en $n+1\geq r+2$ variables $x_{0}, \dots, x_{n}$.  La $k$-vari\'et\'e 
  $X \subset \P^n$  d\'efinie par l'annulation de ces formes  contient un z\'ero-cycle effectif de degr\'e
 $d^{r+1}$.
 \end{lemma}
 \begin{proof}
 Par l'argument donn\'e au lemme \ref{exlemme0} (c), on peut supposer le corps $k$
 infini. Cette hypoth\`ese sera utilis\'ee de fa\c con constante dans ce qui suit.
 Soient  $g_{0},g_{1},\dots,g_{r}$ des formes de degr\'e $d$, \`a coefficients dans $k$,
 dont l'annulation d\'efinit une sous-$k$-vari\'et\'e intersection compl\`ete et lisse $Y \subset \P^n$.
 Soit ${\cal X} \subset \P^n \times \A^1$ le sous-sch\'ema ferm\'e d\'efini par 
 l'id\'eal (homog\`ene en les $x_{i}$) $$(tg_0+(1-t)f_0, \ldots, tg_r+(1-t)f_r)\subset k[t][x_0,\ldots, x_n].$$
 Choisisons des formes lin\'eaires  $L_1,\ldots, L_{n-r-1}\in k[x_0, \ldots, x_n]$ 
 telles que le sous-sch\'ema de
  $\P^n$ d\'efini par l'id\'eal $(g_0,\ldots, g_r, L_1,\ldots, L_{n-r-1})$
  soit fini et \'etale sur $k$.  Alors le sous-sch\'ema ferm\'e
   ${\cal C}'$ de $\cal X$ d\'efini par l'id\'eal
    $(L_1,\ldots, L_{n-r-1})$ est fini et \'etale de degr\'e  $d^{r+1}$ au-dessus d'un voisinage ouvert $U$ of $1$ dans $\A^1$. Soit  ${\cal C} \subset {\cal X}$ l'adh\'erence sch\'ematique de ${\cal C}'\cap \P^n_{U}$ dans  $\cal X$.
   Le sch\'ema  $\cal C$ est propre sur $\A^1_k$ et comme un sous-sch\'ema ouvert dense de $\cal C$ est
  quasi-fini sur $\A^1$, le morphisme $p:{\cal C} \to \A^1$ est fini.
  Comme $p^{-1}(U) \subset {\cal C}$ est \'etale sur $U$, donc r\'eduit,
   son adh\'erence sch\'ematique
  $\cal C$ est aussi  r\'eduite, 
  et chaque composante irr\'eductible de $\cal C$ s'envoie surjectivement sur
 $\A^1$. Ainsi $\cal C$ est plat sur $\A^1$, et donc la fonction
$$
a \to {\rm dim}_{k(a)}{\cal O}_{p^{-1}(a)}
$$
est constante sur $\A^1$. En particulier, le 0-cycle associ\'e au sous-sch\'ema ferm\'e $p^{-1}(0)$ de $X$ a degr\'e $d^{r+1}$ sur $k$.
\qed \end{proof}

\section{Syst\`eme de plusieurs formes, I}

Voici la premi\`ere g\'en\'eralisation du th\'eor\`eme \ref{extheoreme1}.

\begin{theorem}
\label{extheoreme2}
Soient $k$ un corps,  et $f_{0},f_{1},\dots,f_{r}$, avec $r \geq 1$, des formes \`a coefficients dans $k$, de degr\'e $d\geq  1$, en $n+1\geq r+2$
variables. Soit $X \subset \P^n_{k}$ la $k$-vari\'et\'e  d\'efinie par l'annulation de ces $r+1$ formes. 
 Soient $t_{1},\dots,t_{r}$ des
variables ind\'ependantes et $K=k(t_{1},\dots,t_{r})$.
  Soit $W \subset  \P^n_{K}$ la $K$-vari\'et\'e d\'efinie par $f_{1}-t_{1}f_{0}=\cdots=f_{r}-t_{r}f_{0}=0$.
  On a :
  $$I_{red}(X/k)=I_{red}(W/K).$$
   En particulier, la $K$-vari\'et\'e $W$ poss\`ede un z\'ero-cycle de degr\'e 1 si et seulement si la $k$-vari\'et\'e $X$  poss\`ede un z\'ero-cycle de degr\'e 1.
   \end{theorem}
\begin{proof}
Pour les $k$-vari\'et\'es, on omet l'indice $k$. 
Si  $L/k$ est une extension finie de corps  avec $X(L)\neq \emptyset$, alors $W(L(t)) \neq \emptyset$.
Ainsi    l'indice $I(W/K)$  divise $I(X/k)$, et donc l'indice r\'eduit $I_{red}(W/K)$ 
 divise  l'indice r\'eduit $I_{red}(X/k)$.
D'apr\`es le lemme \ref{exlemme0} (c), pour \'etablir l'\'enonc\'e on peut remplacer le corps
$k$ par une extension transcendante pure. On supposera donc le corps $k$ infini.

{\it Supposons d'abord que les formes $f_{i}$ n'ont pas de facteur commun non constant.}

Soient $V=\P^n \setminus X$ et  $I=I(X/k)$. Comme au th\'eor\`eme \ref{extheoreme1}, on a
$\CH_{0}(V)=\Z/I$.
Soit $Z \subset \P^r \times \P^n$ la $k$-vari\'et\'e d\'efinie par
la proportionalit\'e de $(\lambda_{0},\dots,\lambda_{r})$
et de $(f_{0}, \dots, f_{r})$, c'est-\`a-dire par le
syst\`eme d'\'equations $\lambda_{i}.f_{j}-\lambda_{j}.f_{i}=0$ pour $i,j \in \{0,\dots,r\}$.
La fibre de $Z \to    \P^r$ au-dessus du point g\'en\'erique de $ \P^r$ est $W/K$.
Soit $U \subset \P^r \times \P^n$ l'ouvert compl\'ementaire de $Z$.
La projection $q: \P^r \times \P^n \to \P^n$ induit un morphisme $q_{U}: U \to V$.
Soit $q_{Z} :  Z \to \P^n$ le morphisme restriction de $q$ \`a $Z$.
Il induit un isomorphisme $q_{Z}^{-1}(V) \oi V$.
Les fibres de  $q_{U}: U \to V$ au-dessus d'un point $M$ de $V$ sont donc le compl\'ementaire
d'un point
dans un espace projectif $\P^r$.
De fa\c con plus globale, le morphisme $q_{U}: U \to V$  se d\'ecompose comme
$$U \to U_{1} \to V$$
o\`u $q_{1} : U \to U_{1}$ est un fibr\'e en droites affines et $q_{2} : U_{1} \to V$ est un fibr\'e projectif
de dimension relative $r-1$.
Par image directe par morphisme propre on a un isomorphisme $q_{2*} : \CH_{0}(U_{1}) \oi  \CH_{0}(V)=\Z/I$.
Par image inverse par morphisme plat, on a un isomorphisme $q_{1}^* : \CH_{0}(U_{1}) \oi \CH_{1}(U)$.
On a donc un isomorphisme $\CH_{1}(U) \oi  \Z/I$.
On a la suite exacte de localisation
$$ \CH_{1}(Z) \to \CH_{1}( \P^r \times \P^n) \to \CH_{1}(U) \to 0.$$
Par ailleurs les poussettes associ\'ees \`a la  projection $p :  \P^r\times \P^n \to \P^r$ 
et \`a la  projection $q: \P^r \times \P^n \to \P^n$
induisent   un isomorphisme 
$$(p_{*},q_{*}) : \CH_{1}(\P^r \times \P^n) \oi \CH_{1}(\P^r) \oplus \CH_{1}(\P^n) 
= \Z \oplus \Z.$$
Notons $$i : \CH_{1}(Z) \to \CH_{1}( \P^r \times \P^n)  \oi  \CH_{1}(\P^r) \oplus \CH_{1}(\P^n) 
= \Z \oplus \Z$$
l'application compos\'ee.
Sur $Z$, on trouve les 1-cycles suivants.

Un cycle $z_{1}=\P^1 \times R$, o\`u $\P^1 \subset \P^r$ est une droite et
$R$ est un z\'ero-cycle effectif de degr\'e $d^{r+1}$ sur $X$,
dont l'existence est assur\'ee par
le lemme \ref{newlemma}.
L'image par $i$ de $z_{1} $ est $(d^{r+1},0)$.

Soit $s=I(W/K)$. Il existe donc un z\'ero-cycle de degr\'e $s$ sur $W/K$,
et un  tel z\'ero-cycle
 s'\'etend en un $r$-cycle sur $Z$, cycle g\'en\'eriquement fini sur $\P^r$
 de degr\'e relatif $s$. 
 Sa restriction au-dessus d'une droite g\'en\'erale $\P^1 \subset \P^r$
 est un 1-cycle sur $Z$, dont l'image par $i$ est  de la forme $(s,a)$ pour un certain entier $a \in \Z$.
 
 Les formes homog\`enes $(f_{0},\dots,f_{r})$ d\'efinissent un $k$-morphisme
$\sigma : V \to \P^r$ .
Elles d\'efinissent donc une section 
 du morphisme $q_{V} :\P^r_{V} \to V$, restriction de $q$ \`a $\P^r_{V}$,
section dont l'image est dans $Z$ : c'est l'iso\-mor\-phisme inverse
de l'isomorphisme $q_{Z}^{-1}(V) \oi  V$ mentionn\'e plus haut.
L'hypoth\`ese que les $f_{i}$ n'ont pas de diviseur commun non trivial
assure  que la
 $k$-vari\'et\'e $X \subset \P^n$ est de codimension au moins~2. 
Le corps $k$ \'etant  infini, par le m\^eme argument qu'au th\'eor\`eme \ref{extheoreme1},
 on peut donc trouver
une $k$-droite  $L=\P^1 \subset \P^n$  qui ne rencontre pas $X$.
La restriction de $\sigma$ \`a $L$ est donc un morphisme $L \to  Z \subset \P^r \times \P^n$,
dont l'image est un 1-cycle sur $Z$. L'image de ce 1-cycle par $i$
est  $(d,1)$.

Le quotient du groupe $\Z \oplus \Z$ par le groupe engendr\'e par  les trois \'el\'ements $(d^{r+1},0)$, $(s,a)$
et $(d,1)$  est $\Z/J$, avec $J=\pgcd(d^{r+1}, s-ad)$.
De la suite de localisation on conclut  que $I=I(X/k)$ divise $\pgcd(d^{r+1},I(W/K)-ad)$ pour un certain
 entier $a \in \Z$. Ainsi $I_{red}(X/k)$ divise $I(W/K)$.
 Comme par ailleurs $I_{red}(W/K)$ divise $I_{red}(X/k)$, on obtient
 $$I_{red}(X/k)=I_{red}(W/K).$$

 \medskip

{\it Supposons maintenant que $f_{i}=h.g_{i}$ pour tout $i$ avec les   $g_{i}$ homog\`enes de m\^eme degr\'e  
sans facteur commun non constant et $h$ homog\`ene non constant. } 

Soit $X \subset \P^n_{k}$, resp. $X_{1} \subset \P^n_{k}$, resp. $X_{2} \subset \P^n_{k}$ la $k$-vari\'et\'e d\'efinie par l'annulation des $f_{i}$, resp. par l'annulation   des $g_{i}$, resp. par $h=0$.
Soit $W\subset \P^n_{K}$, resp. $W_{1}/K$ la vari\'et\'e d\'efinie par l'annulation des $f_{i}-t_{i}f_{0} \ ( i=1,\dots, r)$,
resp.  par l'annulation des $g_{i}-t_{i}g_{0} \ ( i=1,\dots, r)$.
On a
$$I_{red}(W/K)= \pgcd(I_{red}(W_{1}/K),I_{red}(X_{2,K}/K)) =  \pgcd(I_{red}(W_{1}/K),I_{red}(X_{2}/k)) $$
d'apr\`es le lemme \ref{exlemme0}
et
$$\pgcd(I_{red}(W_{1}/K),I_{red}(X_{2}/k)) = \pgcd(I_{red}(X_{1}/k),I_{red}(X_{2}/k)) =I_{red}(X/k),$$
la premi\`ere \'egalit\'e r\'esultant de $I_{red}(W_{1}/K)=I_{red}(X_{1}/k)$ \'etabli ci-dessus,
 la seconde \'egalit\'e provenant du lemme \ref{exlemme0}.
Ceci ach\`eve la d\'emonstration.
\qed \end{proof}
 
\section{Syst\`eme de plusieurs formes, II}

Voici la seconde g\'en\'eralisation du th\'eor\`eme 1.

\begin{theorem}
\label{extheoreme3}
Soient $k$ un corps, $f_{0},f_{1},\dots,f_{r}$ des formes  non toutes nulles, \`a coefficients dans $k$, 
de degr\'e $d\geq  1$, en $n+1\geq r+2$ variables. Soit $X \subset \P^n_{k}$ la $k$-vari\'et\'e d\'efinie par
 l'annulation de ces formes. Soient $w_{1},\dots,w_{r}$ des
variables ind\'ependantes et $L=k(w_{1},\dots,w_{r})$. Soit $Y\subset \P^n_{L}$
l'hypersurface d\'efinie par  l'\'equation $f_{0}+w_{1}f_{1}+\cdots+w_{r}f_{r}=0$. On a :
$$I_{red}(X/k)=I_{red}(Y/L).$$
 En particulier, la $L$-hypersurface $Y$ poss\`ede un z\'ero-cycle de degr\'e 1 si et seulement si la $k$-vari\'et\'e $X$  poss\`ede un z\'ero-cycle de degr\'e 1.
\end{theorem}
\begin{proof} 
Soit $E_{r}$ l'\'enonc\'e de ce th\'eor\`eme pour $r$ fix\'e et tout corps $k$.
L'\'enonc\'e $E_{1}$ est le th\'eor\`eme \ref{extheoreme1}. Supposons \'etabli $E_{r-1}$.

Supposons $r \geq 2$.
Soient $t_{1},\dots,t_{r} $ des variables ind\'ependantes et $K=k(t_{1},\dots,t_{r})$.
D'apr\`es le th\'eor\`eme \ref{extheoreme2}, on a 
$I_{red}(X/k)= I_{red}(W/K)$,  o\`u la $K$-vari\'et\'e $W \subset \P^n_{K}$ est d\'efinie
par  $$f_{1}-t_{1}f_{0}=\cdots=f_{r}-t_{r}f_{0}=0.$$
Soient $s_{2},\dots,s_{r}$ des variables ind\'ependantes et $F=K(s_{2},\dots,s_{r})$.
D'apr\`es $E_{r-1}$, l'indice r\'eduit de $W$ sur $K=k(t_{1},\dots,t_{r})$ est \'egal \`a l'indice
r\'eduit sur $F$ de l'hypersurface $T$ d\'efinie dans $\P^n_{F}$ par
$$ (f_{1}-t_{1}f_{0}) + s_{2} (f_{2}-t_{2}f_{0})+ \cdots + s_{r }(f_{r}-t_{r}f_{0})=0.$$
Ceci se r\'e\'ecrit
$$f_{1} - (t_{1}+s_{2}t_{2}+\cdots+s_{r}t_{r})f_{0} +s_{2}f_{2}+ \cdots +s_{r}f_{r}=0.$$
Soient $w_{1}=-1/(t_{1}+s_{2}t_{2}+ \cdots +s_{r}t_{r})$ et, pour $i \geq 2$,
$w_{i}=- s_{i}/(t_{1}+s_{2}t_{2}+ \cdots +s_{r}t_{r})$.
\vfill\eject
L'\'equation de l'hypersurface $T \subset \P^n_{F}$  s'\'ecrit alors
$$f_{0}+w_{1}f_{1}+ \cdots +w_{r}f_{r}=0.$$
L'inclusion 
$$k(w_{1},\dots,w_{r}, t_{2}, \dots, t_{r}  ) \subset k(t_{1},t_{2}, \dots, t_{r},s_{2},\dots,s_{r})=F$$
est   une \'egalit\'e. L'extension $F=L(t_{2}, \dots, t_{r} )$ est transcendante pure.
D'apr\`es le lemme \ref{exlemme0}, l'indice r\'eduit sur $L$ de l'hypersurface
d\'efinie par $$f_{0}+w_{1}f_{1}+ \cdots +w_{r}f_{r}=0$$ dans $\P^n_{L}$
est \'egal \`a l'indice r\'eduit de cette hypersurface sur $F$.
Ceci ach\`eve la d\'emonstration.
\qed \end{proof}

\begin{remark}
A. Pfister, J.W.S. Cassels et D. F. Coray  (voir les r\'ef\'erences dans \cite{C}) ont donn\'e
donn\'e des exemples d'inter\-sec\-tions compl\`etes de trois
quadriques $f_{0}=f_{1}=f_{2}=0$ dans $\P^n_{k}$ 
(sur un corps $k$ de caract\'eristique diff\'erente de 2)
qui poss\`edent un z\'ero-cycle de degr\'e 1 sans poss\'eder de point rationnel.
La quadrique  $f_{0}+t_{1}f_{1}+t_{2}f_{2}=0$ sur le corps $k(t_{1},t_{2})$ poss\`ede 
alors un z\'ero-cycle de degr\'e 1. Comme c'est une quadrique, un th\'eor\`eme de
Springer \cite{S} assure que cette quadrique admet un point $k(t_{1},t_{2})$-rationnel.

On voit ainsi  que le   th\'eor\`eme \ref{extheoreme3}
ne vaut pas lorsque l'on remplace les z\'ero-cycles 
de degr\'e~1 par des points rationnels : le th\'eor\`eme
d'Amer et Brumer ne s'\'etend pas \`a un syst\`eme de 3 formes.
\end{remark}

\begin{acknowledgement}
Les auteurs savent  gr\'e au rapporteur
de ses lectures attentives.
Marc Levine remercie la  NSF (grant number DMS-0801220) et
la fondation Alexander von Humboldt. \end{acknowledgement}

\end{document}